\documentclass[a4paper,12pt,leqno]{article}

\usepackage{geometry} 

\usepackage{tikz-cd}\usepackage{bussproofs}
\usepackage{amssymb,amsmath}
\usepackage{proof}
\usepackage{psvectorian}
\usepackage{comment}
\usepackage{CJKutf8}
\usepackage[utf8]{inputenc}
\usepackage{graphicx}
\usepackage{graphicx}
\graphicspath{ {images/} }
\usepackage{qtree}
\usepackage{mathtools}
\usepackage{hyperref}
\usepackage{listings}
\usepackage{indentfirst}
\usepackage{titlesec}
\usepackage{yhmath}

\DeclareUnicodeCharacter{00AC}{$\neg$}

\DeclareUnicodeCharacter{03B5}{$\epsilon$}

\DeclareUnicodeCharacter{2200}{$\forall$}

\DeclareUnicodeCharacter{2203}{$\exists$}

\DeclareUnicodeCharacter{2229}{$\cap$}

\DeclareUnicodeCharacter{222A}{$\cup$}

\DeclareUnicodeCharacter{2282}{$\subset$}

\DeclareUnicodeCharacter{03C9}{$\omega$}

\DeclareUnicodeCharacter{2218}{$\circ$}

\DeclareUnicodeCharacter{207B}{$^{-1}$}
\numberwithin{equation}{section}
\newtheorem{defn}[equation]{Definition}

\newtheorem{rem}[equation]{Remark}
\newtheorem{exm}[equation]{Example}

\newtheorem{notat}[equation]{Notation}
\newtheorem{newpar}[equation]{}
\newtheorem{xdefn}{Definition.}
\newtheorem{xproposition}{Proposition.}
\newtheorem{xcorollary}{Corollary.}
\newtheorem{xrem}{Remark.}
\newtheorem{xexm}{Example.}
\newtheorem{xlemma}{Lemma.}
\newtheorem{xtheorem}{Theorem.}
\newtheorem{xnotat}{Notation.}
\newtheorem{xnewpar}{\it}
\newtheorem{xproof}{{\it Proof. }}
\newtheorem{xproofof}{{\it Proof}}

\newenvironment{remark}{\begin{rem}\em}{\end{rem}}

\newenvironment{newparagraph*}[1]{\begin{xnewpar}\hspace*{-1.5mm}{#1}. \rm}{\end{xnewpar}}
\newenvironment{definition*}{\begin{xdefn}\em}{\end{xdefn}}
\newenvironment{remark*}{\begin{xrem}\em}{\end{xrem}}
\newenvironment{example*}{\begin{xexm}\em}{\end{xexm}}
\newenvironment{notation*}{\begin{xnotat}\em}{\end{xnotat}}
\newenvironment{proposition*}{\begin{xproposition}}{\end{xproposition}}
\newenvironment{corollary*}{\begin{xcorollary}}{\end{xcorollary}}
\newenvironment{lemma*}{\begin{xlemma}}{\end{xlemma}}
\newenvironment{theorem*}{\begin{xtheorem}}{\end{xtheorem}}

\titleformat*{\section}{\large\bfseries}

\begin{document}
\author{Clarence Lewis Protin}
	
	\title{Modern Definition and Ancient Definition (draft)}
	\maketitle 
	
	\begin{abstract}
	In this essay we examine some aspects of the classical theory of definition as codified in Aristotle's \emph{Topics} and Porphyry's \emph{Eisagogê} in the light of the way definition is carried out in modern mathematical practice. Our goal is to contribute to the understanding of the alleged gap existing between ancient and modern logic and science as well as the reasons behind allegations of inadequacy and lack of sophistication in the ancient theory of definition. Also to investigate the possibility of a co-interpretation between modern mathematical definitional practice and ancient definitional practice in particular in the light of topos theory. We find the ancient definitional practice asks relevant and overlooked questions about modern mathematical practice which apparently have escaped current philosophical and mathematical logical literature. We also present some general considerations about the structure and development of theories as these relate to the theory of definition. 
	
	\end{abstract}

\section{Introduction}
	
This paper aims to contribute to the historical and philosophical study of the theory and practice of definition by presenting a formal mathematical outline of the ancient Aristotelic theory and then studying the practice of definition in modern mathematics in light of the ancient theory while also highlighting its divergences and contrasts. 	
This paper is organized as follows. In section 2 we present a general formalization of the main concepts of the ancient theory of definition such as found in Aristotle's \emph{Topics} and Porphyry's Eisagogê and point out some obscure aspects and open questions. In section 3 we present a more sophisticated mathematical model for a fragment of the ancient theory.
In section 4 we investigate the practice of definition in modern mathematics in the light of the ancient theory and study the persistence and covert presence of the classical concepts as well as striking departures from them. We show how the modern and ancient practices ask important questions about each other. 
In section 5 we bring together section 2, 3 and 4 to suggest how the ancient and modern theories might be seen as unified in the context of modern topos theory.
In section 6 we look at the larger picture of how definitions in mathematics and science articulate with each other to form theories and some of logical and epistemological questions are raised which can also be traced back to fundamental preoccupations of ancient philosophy. 

\section{Outline of some aspects of the ancient theory of definition}

In order to effect a comparison between the ancient and modern practices of definition is it convenient to give a succinct summary of the ancient practice. The oldest systematic source which has come down to us is Aristotle's \emph{Topics} which in a very rich and complex text which encompasses grammatical, semantic, pragmatic, epistemic and ontological concerns as much as logical ones in the modern sense (see \cite{slo, prim, pro3} for a justification of the reading of the \emph{Topics} we follow here as well as general references\footnote{all references to the \emph{Topics} follow Brunschwig's edition \cite{topiques1,topiques2}.}). There is also a much later, shorter and perhaps somewhat syncretic  version of this Aristotelic theory presented in Porphyry's \emph{Eisagogê}\cite{por}. In this section we outline some basic important aspects and concepts of the ancient theory as presented in these two texts. Our presentation will be a formal mathematical one.

 For some arguments for the legitimacy of the formalization of historical philosophical systems (what is called the hermeneutical  'virtuous circle') see for instance \cite{lam}. We remark that for the purpose of comparison with the practice of definition in modern mathematics having such a formal version of the ancient theory can only be an advantage, as it provides a quick and accessible introduction to the ancient theory for the modern mathematician or philosopher. Also, having a formal version sheds light on the obscurities and difficulties of the classical theory which cannot afford to overlooked if we wish to make an adequate comparison and  co-interpret with the modern theory.



  We start with non-empty set $U$ with a disjoint decomposition $U = X \cup S$ where $X$ and $S$ are non-empty and $S$ is finite. The elements of $S$ are called \emph{terms}\footnote{For a good discussion of Aristotle's notion of term \emph{oros} see \cite{ham}. Our formalization is meant to be general enough to be interpretable both syntactically or semantically.} and those of $X$ \emph{substance individuals}. $X$ may be countably infinite. There is a fundamental binary relation on $U$ which is denoted by $\overline{\in}$. It is required to be transitive\footnote{this transitivity is the ancestor of the Barbara syllogism of the \emph{Analytics}.} and antireflexive
$\forall u \in U(u\overline{\notin} u)$. The elements of $X$ are atoms, that is, for $x \in X$ there is no $u \in U$ such that $u\overline{\in} x$. Also for every element of $x \in X$ there is some $s \in S$ such that $x\overline{\in} s$. The following definition is fundamental:

\[ u\in w :\equiv u\overline{\in} w\,\&\, \sim\exists s\in S(u\overline{\in}s\overline{\in} w)\]

If $u \in w$ we say $w$ is a \emph{proximate genus} of $u$ (it will follow from other axioms that it must be unique) and $u$ is a \emph{proximate species} of $w$ except in the case in which $u \in X$ where we say that $w$ is the \emph{infima species} of $u$.

The relation $\overline{\in}$ of subordination between species and genus is fundamental. One can argue that Book IV of the \emph{Topics} is dedicated to the topics which follow from its basic properties.

Topics are a series of sentences which when instantiated take the form of an implication, the consequent being either the affirmation of negation of the an instantiation of the genus, property or definition relation. The expression of the transitivity of $\overline{\in}$ is itself an example of a constructive (affirmative) topic.

 Transitivity is postulated in the justification of the very first topic 120b15-25 (and from several others such as 122a30-35). Antireflexivity follows, for instance, from a topic in 121b10-15: equal terms have the same extension and therefore cannot be subordinate to each other, for Aristotle postulates that if $u\overline{\in} v$ then the extension of $v$ must be larger than $u$ (this will be formalized as an axiom further ahead).  The last axiom which in itself is evident in the context of Aristotelic philosophy. The postulate of the finitude of $S$ follows from the empirical fact of the finitude of linguistic expressions. It is also implied in many topics in Book IV, for example 122a3-6\footnote{of course for a general $S$ we could account for these topics by postulating a Noetherian condition on ascending chains.}. 

A very important axiom (T), found in the last topic in 121b24-30 gives $\overline{\in}$ its tree-like structure:

\[\text{if $u \overline{\in} v$ and $u \overline{\in} w$ then $v\overline{\in} w$ or $w\overline{\in} v$} \tag{T}\]

This tree-like structure becomes even more definite in the  topic 123a30-32 giving the important principle that \emph{a genus must have more than one proximate species}.

Equipped with this axiom we can now introduce the 'categories' as follows:

\begin{itemize}
\item There are ten maximal elements $\top_1,...\top_{10}$ in $U$.
\item There is one and only one maximal element $\top_1$ for which $x \overline{\in} \top_1$ for all $x \in X$.

\end{itemize}

Regarding the maximal elements it interesting to note that 'universal' terms like 'One' and 'Being' are in 127a25-25 excluded by Aristotle from being genera. Thus $S$ may contain isolated elements which have no species and which belong to no genera. It is patent that the relation $\overline{\in}$ is only a particular case of the general Aristotelic concept of 'predication'.

If follows from the above axioms that there is a decomposition of $U$ into ten disjoint non-empty subsets $U_i$ for which $X \subset U_1$ and if $u,v$ belong to different $U_i$'s then we cannot have $u\overline{\in} v$ or $v\overline{\in} u$ (this is the basis of many topics in Book IV: 120b36-39, etc.). The ten maximal elements corresponding to the ten $U_i$ are Aristotle's famous ten categories. It also follows from the axioms that for each $x \in X$ there is one and only one $s \in U_1$ such that $x \in s$. 

\begin{remark} How do we interpret the atoms for $U_i$ with $i\neq 1$ ? Does it make sense to speak of individuals or infima species for categories other than substance ?
\end{remark}

\begin{remark} Do all chains $s_1 \in s_2 \in...\in \top_i$ with $s_1$ atomic have the same length ? If so we can divide each $U_i$ into a series of levels $U_i^1,...U_i^{L_i}$. What can we say in general about the cardinalities of the $U_i^j$ as $j$ increases ? This is important for studying in more detail the structure of \emph{Porphyry's tree}.
\end{remark}

We now attempt to formalize \emph{difference}. To this end we postulate a partial function

\[\otimes: U\times S \rightharpoonup U\]

which represents the way a term acts on a genus as a difference to yield a species. It can be thought of as a partial $S$-module multiplication acting on each $U_i$. This operation must satisfy:

\[\text{$u \in U_i, v \in U$ implies $uv \in U_i$ whenever defined} \]

where we denote $\otimes$  henceforth by concatenation.
We also require that $xv$ is never defined for $x \in X$.

  Given $s \in S$ we denote by $\Delta(s) \subset S$ the set of all $u \in S$ for which $su$ is defined. Thus to each genus $s$ there is a set of terms which can be used to carve up $s$ into subordinate species.

 Implicit in the entire ancient theory of definition is the possibility of giving a definition to each term, hence:

\begin{itemize}
	
\item For all $s,t \in S$ such that $s\in t$ there is a $d \in S$ such that  $s = td$.
\end{itemize}

We call $d$ the difference of $s$ relative to $t$. The expression $td$ is itself the \emph{definition} of $s$.

Many topics in Book IV (122b18-20, etc.) are based on the postulate that 

\[ s = td \rightarrow s\overline{\notin} d\]

The differences of a species (relative to some genus) are radically distinct from the genera to which it is subordinate.

\begin{remark} Is $d$ unique ? Is difference unique in a correct definition for Aristotle and Porphyry ? This is suggested for instance by 141a 35-36. And which category does $d$ belong to ? For $s,t$ in the category of substance it seems that $d$ can be for instance a quality\footnote{144a18-19: \emph{dokei d'hê diaphora poion ti sêmainein.}} but never a substance as well (cf. 143a29-33). But what if $s,t$ are not in the category of substance, what should difference be ?  Aquinas wrote regarding this in \emph{De Ente e Essentia}, VII: \emph{Definitionem autem habent incompletam, quia non possunt definiri nisi ponatur subjectum in eorum definitione}.
\end{remark}	

It is difficult to formalize exactly what was meant by the \emph{extension} of a term in the ancient theory. In the \emph{Topics} sometimes Aristotle seems to mean the set of individual substances of which the term is predicated of and sometimes the set of terms of which the term is predicated of.  A tentative approach would be to postulate an extension function:

\[ E:  S \rightarrow PX \]

where $PX$ is the set of subsets of $X$. The function $E$  must satisfy the condition discussed earlier: $s \overline{\in} t$ implies that $Es \subsetneq Et$\footnote{this condition permeates the Topics. It is used in particular in 127a25-35 to prove that 'One' and 'Being' cannot be genera.}.

The interpretation of extension for the category of substance is clear: the extension of 'swan' is the set of all swans. But what is the extension of terms belonging to the categories of quality, quantity or relation ? The extension of 'white' would have likely been considered the set of white  substances\footnote{In 127a20-25 Aristotle argues that \emph{to leukon} is not a genus because its subordinate individuals do not differ in species. It is not clear if we should take this term as 'white' or 'white things' and thus it is not clear if white things are to be taken as the extension of 'white'. }.

\begin{remark} What would the extension of a \emph{relation} be, ? Let $U_r$ constitute the terms in the category of relation. Then in this case we must emend our definition to 
	\[ E' : U_r \rightarrow P(X\times X)  \]
	For instance the extension of the relation of marriage would the  of pairs $(a,b)$ in which the individuals $a$ and $b$ are married. We can also ask how definitions were carried out for relations. Consider the relation of fatherhood. How do we define the relation of being a paternal grandfather without using multiple generality (see \cite{bob1} for a discussion of multiple generality in Stoic logic) ? Also what would be the genus and difference here ? It is difficult to see how our partial multiplication function could express this without further logico-algebraic refinement.
\end{remark}

The previous remark introduces one of the most outstanding problems in the interpretation and formalization of the ancient theory of definition. The problem is that \emph{differences may need to be grammatically and logically complex terms}. It seems we need to further refine out set $S$ to make the distinction between simple and complex terms and to introduce various logical and grammatical operations to generate complex terms. It would be interesting to make a systematic study of the logical and grammatical forms of difference terms in the \emph{Topics}. In the present section we consider only a basic form of term combination\footnote{Aristotle considers the definitions of complex (\emph{sumpeplegmenos}) terms in the \emph{Topics}, for instance in 148b23. The rules he lays down are algebraically suggestive. The first one might be translated as a cancellation law:  $sd_1 =  td_2d \rightarrow d_1 = d$ where $d_1 = d$ must be a definition.}.

Let us consider the challenge of formalizing the two remaining \emph{predicabilia}, property (\emph{idion}) and accident. This is a complex subject and there are many nuances (i.e. the difference between separable and inseparable accident) and apparent inconsistencies relating to these terms in the \emph{Topics}. Porphyry offers a simpler and more intelligible account, specially in function of his methodic comparison of these terms with genus and difference. Here are tentative definitions of property and accident.

We say that $t$ is a \emph{property} of $s$ if $Et = Es$ and $t$ cannot be written in the form $t = a d$ for some $a$.

If $Es \subsetneq Et$ and it is not the case that $s\overline{\in} t$ then $t$ is an (inseparable) \emph{accident} of $s$.

\begin{remark} We must investigate how a property differs from difference. Indeed we have not postulated anything about the extension of differences. We could say, of course, that property is not a difference because it will never be part of a decomposition $s = u \otimes t$. But Porphyry states that difference is different from property because such a term can be the difference of various species while a property is only the property of one species. Also property and its species have a kind of (intensional) convertability which difference lacks. 
\end{remark}

In order to formalize the \emph{comprehension} of a term we need to introduce a conjunction partial function

\[ \times : S \times S \rightharpoonup S\]

which satisfies the properties:

\begin{itemize}
\item  $t,s \in U_i$ implies that $t\times s \in U_i$.
\item  $t (s \times u) = (t s) u$\footnote{Note that this is similar to the definition of an action or module in algebra.}.
\item  if $t \in s \in u$ and $s = ut_1$ and $t = st_2$ then $t = u(t_1\times t_2)$.
\end{itemize}

Syntactically this operation corresponds to concatenation of differences, for instance 'featherless biped'. Thus we must consider $S$ as representing both simple and compound expressions (but of bounded length, maintaining the finitude of $S$: we can realistically consider that $\times$ is simply undefined for large compound expressions).

 Since $S$ is finite there will be simple or indecomposable terms and we can define in the expected way the maximal decomposition of a term. 
 
 The reason for introducing $\times$ is that now given a chain  \[ s_1 \in s_2 \in...\in s_k\in...\in \top_j (= s_n)\]
 with $s _i = s_{i+1}d_i$ we can define the \emph{comprehension} $c$ of $s_1$ as $c = d_1 \times d_2 \times...\times d_{n-1}$ and we will have that $s_1 = \top_j c$. In Porphyry's example the difference of man relative to the supreme genus would decompose into mortal, rational and sensible - all of which would be components of the comprehension of the term man.
 
 \begin{remark} Note that it follows from the above postulates that we can have equalities $s = td$ which are not, strictly speaking, definitions when it is the case that $s \overline{\in} t$ but not $s \in t$.  This is in accordance with Aristotle's rule in 143a19-20 that $s$ must be placed in the proximate genus \emph{eis to eggutatô genos theinai}.
 \end{remark}

This is only meant to be a brief sketch of the formalization. There are many important developments and details which remain to be discussed.

 The \emph{Topics} offer a plethora of refinements and details with regards to the general theory which are not discussed in Porphyry's purposely brief and introductory text.  Some interesting aspects are the required non-redundancy of $m$ for a definition $t = s  m$. For instance $m$ cannot admit a decomposition $m = m_1 \times m_2$ in which $m_1$ has the same extension as $t$. Also repetitions of factors cannot occur in a $\times$-decomposition of $m$ (we shall return to this in detail in section 4). And if we were to introduce term-forming operations which can instantiate relations then there is the problem of circularity in definitions which is again an important aspect of the \emph{Topics} which we will return to in section 4. And there is the theory of opposition: if $t = s m$ and $t \in s$  there can be a situation of dichotomy in which the other species $t' \in s$ is given by a kind of complementation operator $t' = s m^\circ$. There are many other refinements and nunances considered by Aristotle.


We note that the present formalization allows two interpretations of universal quantification. One involving the containment of extensions and the other an 'intensional' interpretation in terms of $\overline{\in}$. In the following section we give a more sophisticated extensionalist model.

\section{Another model for the ancient theory of definition}

In the following model we consider only the category of substance (and touch upon the relations involving individual substances). This model is an extensionalist interpretation\footnote{But we shall see later how this aspect might be overcome.} and is based on a bottom-up approach. The idea is to start with a set $X$ of all substance-individuals (from now on referred to just as individuals) and successively construct the infima species and higher and higher genera by the operation of collecting elements of a set into a finite partition\footnote{Recall that given a set $Y$ a partition of $Y$ consists of a set $P$ whose elements are subsets $S_i \subset Y$ for $i=1,..,n$ such that $S_i \cap S_j \neq 0$ for $i\neq j$ and $\bigcup_{i=1,...,n} S_i = Y$}. 
The basic idea is as follows. We start with the set $X$ of individuals. Then we note that the extensions of all infima species determine a partition $S_1$ of $X$. All elements of $x$ are in one and only one infima species and each infima species $\sigma$ determines a non-empty subset of $X$. Now $S_1$ consists of subsets of $X$,  it is thus itself a subset $S_1 \subset PX$ where $PX$ denotes the set of subsets of $X$. But since $S_1$ is a subset of $PX$ it is an element of the set of all subsets of $PX$, that is, $PPX$. So we have $S_1 \in PPX$. Note that a partition of $X$ is thus an element of $PPX$, a set of subsets of $X$. But it is easy to see that not all elements of $PPX$ are partitions. The partitions of $X$ are a proper subset of $PPX$ (we also use the notation $P^2X$) which we denote by $G X$. Thus $S_1 \in G X$. 
Now we also have $S_1 \subset PX$ thus $S_1$ is a set and we can repeat the whole construction by considering a partition $S_2$ of the set $S_1$, this to capture the proximate genera of the the set of infima species. Thus we can consider a $S_2 \subset PS_1 \subset PPX$ for which $S_2 \in G S_1 \subset PPS_1$. We can repeat this construction to obtain a sequence $S_1, S_2,...,S_i,...$ for which

\[  S_i \subset PS_{i-1} \subset P^i X \text{  and  } S_i \in GS_{i-1} \subset P^2 S_{i-1}\]

This gives us a model of the universe of all species and genera in the category of substance\footnote{this should be compared to the suggestion of the existence of levels in remark 2.2.}. Given a $\sigma \in S_i$ its immediate species consists of its elements when viewed as a set and its immediate genus consists of the element of $S_{i+1}$ to which it belongs (this being a partition of $S_i$). Up until now we have used the following postulates from the ancient theory of definition: that the partitions are finite (But $X$ need not be finite); that the elements of the partition are non-empty subsets. It is easy to see that given our previous construction, starting from a $x \in X$ there is a uniquely sequence \[x \in \sigma_1 \in \sigma_2 \in \sigma_3 \in ...\in \sigma_i \in \sigma_{i +1} \in ...\]
such that $\sigma_{i+1} \in S_{i+1}$ is the proximate genus of $\sigma_{i}$. We call the above sequence the \emph{Porphyry sequence } of x.

We are, of course, missing a fundamental postulate of the ancient theory: that this sequence must stabilize  at a supreme genus $\top$. That is, there must be a $n$ for which $S_n$ is the trivial partition  $\{S_{n-1}\}$. With this postulate we can verify that our model captures all the fundamental properties of Porphyry's account of the genus-species relation for the category of substance.

We view our constructed hierarchy as expressing the repeating operation of \emph{synthesis}, of bringing together the many into wholes based on common attributes, which is repeated from the ground up in an ascending hierarchy. The opposite operation (corresponding to analysis) is that of taking the extension of a genus. That is, to each $\sigma \in S_i$ we can associate a subset $E\sigma \subset X$ such that for any $x \in E\sigma$ its Porphyry sequence involves $\sigma$ (it may continue). Given any set $Y$ there is a well-defined map 
\[ \epsilon_Y : PPY \rightarrow PY\]
given by $\epsilon_Y (S) = \bigcup_{s \in S} s$ for $S \in PPY$. 
Now consider $S_i \subset PP(P^{i-2} X)$ and suppose that $\sigma \in S_i$ for $i\geq 3$. Then we define

\[ E\sigma = ... \epsilon_{P^{i-3} X} \circ \epsilon_{P^{i-2} X}(\sigma)  \in PX\]

that is, we repeatedly apply $\epsilon$\footnote{Where no confusion arises we omit the subscript.} until arriving at an element of $PX$, i.e., a subset of $X$ which will be the extension of $\sigma$. This procedure justifies the description of the taking of an extension as an analysis, it has a katabolic nature as $\epsilon$  is destroying information.

But how do we model difference ? Given a $\sigma \in S_i$ we have that $\sigma \in PS_{i-1} \subset PPS_{i-2}$ so $\sigma \in PPS_{i-2}$. Thus if $\sigma' \in \sigma$ is a species then $\sigma' \in PS_{i-2}$ so that $\sigma' \subset \epsilon \sigma$. We thus propose that the difference of $\sigma'$ relative to its proximate genus $\sigma$ is given by the characteristic function 
\[ \phi_{\sigma'} : \epsilon \sigma \rightarrow \{0,1\}\]
Recall that there is a bijective correspondence between subsets $W \subset Y$ and functions $f_W: Y \rightarrow \{0,1\}$ for which $f(x) = 1$ iff $ x \in W$. Further investigation of this model as well as the modeling of property and accident will be theme for future work.

If we repeated the constructions in this section starting not from $X$ but from the Cartesian product $X\times X$ then we might obtain an adequate model for a particular kind of relation, that of relations between substances.

\section{Ancient and Modern Definition}

There can be no doubt that \emph{definitions} played a key role in ancient philosophy and ancient philosophical debate.  A common view is that such practice rested upon assumptions which would be seen as debatable to us today, such as the authority of language\footnote{In the \emph{Topics} there appears the implicit belief that common language holds within its combinatoric capacities (exercised in debate) the possibility of yielding valid knowledge of things. } or the ability of the mind to non-discursively and non-sensuously grasp the essence of things. It is held that the ancient process of seeking definitions is irrelevant to scientific progress\cite{kne}.  Modern predicate logic seems at once to be simpler and more powerful (and less ontologically committed) than the clumbsy natural-language based mechanics of ancient logic\footnote{though it also collapses when trying to deal with changing contingent  perceptible singulars, with the sense of proper names, the domain of \emph{generation and corruption} which the ancient Greeks instinctively avoided. }.  

And yet the basic structure of the ancient theory of definition seems to be suprisingly persistent in the history of western thought. Consider taxonomy and cladistics in biology and biochemistry or the structures used in object oriented programming,  knowledge representation and lexicology. The template of the ancient theory of definition permeates much of modern scientific knowledge. But the focus of the present paper is on mathematics.

  To define something in mathematics\footnote{We shall assume in what follows that we are in the framework of standard set-theoretic mathematics, ZFC or some equivalent system. At the end of the section we touch briefly on the type-theoretic framework.} we start out with a previously given set (its 'matter', so to speak) and apply a condition to this matter to obtain  a 'differentiated' proper subset of this set, characterized uniquely by this condition within the previous enveloping set. Thus groups form a proper subset of monoids given by the condition that there exists a unique inverse for each element. Note that this condition only 'makes sense' for monoids.    Let $Mon$ denote the class or set of Monoids and $Grp$ the set or class of groups. Let $\phi(x)$ be the condition on a monoid $x$ of having all elements possessing a unique inverse.  Then we can write

\[ Grp = \{ x \in Mon : \phi(x)\}\]

which we can read \emph{a group is a $\phi$-satisfying monoid}. 
In the notation of section 2 we could write $ Grp \in Mon$ and $Grp = Mon\otimes \phi$. Which is to say: $Grp$ is a species of $Mon$ with difference $\phi$.

But what criteria are used in the choice of $\phi$ ? Clearly there are infinitely many equivalent choices.  Aristotle gives a number of conditions in \emph{Topics} Book VI. One condition is non-redundancy, given in 140a33-35. In modern terms this means that $\phi$ cannot be of the form $\phi_1 \,\&\, \phi_2$ in which $\phi_1$ alone is already a definition of our object (and in particular co-extensional). Aristotle's non-redundancy condition could be interpreted in a more general way, for instance as rejection redundancy of the form $\phi_1 \vee \phi_2$ in which $\phi_1$ is a definition; we can also think of redundancy situations in which the guilty subformula appears within a quantifier. Thus Aristotle's condition challenges us to find precise formal account of the non-redundancy of a defining formula. Then there is the fundamental condition of impredicativity or non-circularity posited in 142a34: the \emph{definiendum} cannot occur in the \emph{definiens}. We will return later to the question of impredicativity later when we discuss type theory and a case of impredicative definition found in Sextus Empiricus.

   Let $Ab$ be the class of Abelian groups. Then $Ab = Grp\otimes \psi$ where $\psi$ is the commutativity condition. It would not be a good definition to define $Ab$ directly in terms of $Mon$.  Thus to logical non-redundancy and non-circularity we must add a certain minimality condition which likewise finds its place in the \emph{Topics} and which corresponds to the \emph{proximality} condition of $\in$ in terms of $\overline{\in}$.

But what guarantees that $Grp$ is the proximate genus of $Ab$ ? What prevents us from considering the proximate genus of $Ab$ to be Hamiltonian groups $Ham$ (those in which all subgroups are normal) ?  And given any non-singleton set or class what guarantees that it is an infima species ?  For instance for any class of groups we might always think up some further division of this class.   \emph{This suggests that, breaking with classical dogma,  we should consider   the concepts of immediate species and proximate genus as relative and provisional reflecting the dynamic, growing, self-correcting nature of human knowledge}. But before all else a crucial difference is that in modern mathematics there is no finiteness condition on the cardinality of terms as postulated in the ancient theory (see section 2).

If we look at 'difference' in mathematical definition then we see that it often seems to function Platonically or 'dichotomically': groups are divided in abelian and non-abelian groups and so forth.  In general we would expect that each genus be decomposed into a finite (again a condition in the \emph{Topics}) set of disjoint species as in the classification of finite simple groups. Mathematicians call such decompositions of a given genus a 'classification'.  But note that it can be argued that in each species of simple groups there are still infinitely different 'kinds' of groups, for instance cyclic groups of prime order. Thus either we consider these as 'individuals' or we must break with the ancient theory.

  Mathematicians also find 'characterisations' of given classes of structures.  What distinguishes characterisations from definitions ? Characterisations are  the modern descendants of the Aristotelic 'property' or \emph{idion}. While modern mathematics certainly rejects a dogma of the uniqueness of the definition (a problem we raised in Remark 2.3), there is still the philosophical important question of what constitutes a good definition and what should be relegated to being considered a characterisation.
  
  In the \emph{Topics} (as seen in our account in section 2) it appears as if differences where somehow organically connected to genera, as if they were innate to their genera. One cannot just take any condition our of the blue and carve up a genus - and this is said to be one of the main difference of Aristotle's theory of definition from that of Plato.  In this respect we observe that a collection of conditions applicable to a given structure $M$ does not in general decompose $M$ neatly into a tree-like structure (strictly speaking we would expect something isomorphic to a Boolean lattice) of genera and species.  It is interesting to look at figure 5 on p.23 of Steen and Seebach's \emph{Counterexamples in Topology}\cite{steen}  for a sobering example of the messy situation for compactness conditions on a topological space.

  We can understand the concept of group without the concept of commutativity but we cannot understand the concept of commutativity without that of a magma, that is, without some genus of 'group'. For how do we formalise $X$ is commutative, $Com(X)$ ? Clearly this will involve $Magm(X)$. This at least suggests that commutativity is a legitimate embodiment of difference.

  In mathematics not every class proposed as a proximate genus of a another class is legitimate  or 'natural'. As in the \emph{Topics}, the immediate species must have some structural similarity and share some relation between each other\footnote{for instance, form a category or more specifically a category with basic properties such as admitting (finite) products or coproducts.}. The artificial class of finite algebraic structures with an odd number of elements is not a good genus.
  Thus the ancient theory poses the question: why is it that most mathematical structures happen to form categories ?
  
  In the classical theory of definition there is, as we saw, the concept of \emph{infima species} beneath which there are only individuals which cannot differ in species.  What would the counterparts in modern mathematics be ? Should objects like $\mathbb{Z}$ or $\mathbb{R}$ be considered like individuals or infima species ? If infima species then their definition would seem to have to be given in terms of the categoricity of their theories (in the sense of modern model theory). But what if they are defined only 'up to isomorphism' ? 
  And importantly these objects would seem to violate axiom (T)  by being members of non-subordinate distinct genera: they are both concrete examples of a group and a ring.
  And what gives such structures their 'individuality' or infima species quality: is this, in model-theoretic terms their 'categoricity' ? Are mathematical structures defined only 'up to isomorphism', thus ruling out individuality ?  And how do we interpret the set-theoretic elements of these structures ? Are they like the parts of individual substances ?

  Aristotle's view that there is no knowledge of the individual is well known.
  Indeed  concrete space-time individuals continue to elude us. Not only cannot we define them but we cannot, it seems, uniquely identify them. Consider the problem of giving an account of proper names. Ultimately every  unique identifier is only so relative to a certain historical community at a certain place and time. The community has to assumed as given. There are also no absolute coordinates in cosmology which would allow us to identify a given moment of time or location in space. Coordinates are always relative to a certain community and system of knowledge or belief. All these reference systems are incomplete and themselves lacking foundation. Pragmatics is simply the logic of concrete existence. So in this respect modern research only confirms the ancient position.
  
  What about the modern counterpart of 'accident' ?  According to Porphyry the  difference between 'property' and 'inseparable accident' is that the inseparable accident is applicable to other species. Perhaps the modern counterpart of accident is simply any lemma about a given mathematical structure. For instance: Artinian rings are Noetherian rings: so being a Noetherian ring is an accident of Artinian ring. This example is also interesting because the definition of Artinian and Noetherian appear to be dual to each other and \emph{intensionally} distinct, though extensionally one class is contained in another.

The concept of 'opposition' plays an important role in the \emph{Topics}.
Modern negation  is only one type of opposition. The modern counterpart of ancient opposition - which is doubtlessly rooted in the structure of natural language semantic spaces - is to be sought in more subtle kinds of dualities and symmetries which abound in mathematical structures. This is topic for further research.

In Aristotle there is a basic division between  \emph{essential} predication, \emph{kat'hupokeimenou} (that which applies only to individual substances and tells us what they are) and \emph{en hupokeimenô(i)} predication, relating to that which resides within individual substances. We call this last attribute predication. It seems that there are three degrees of attribute predication which differ in their degree of essentiality, on how close the come to essential predication (without ever coinciding with it).
An example of essential attribute of an individual in the \emph{Topics} is 'featherless biped' or 'rational'.  There are also permanent but non-essential attribute like the ability of laugh. In Book V of the Topics the non-essentiality of property is unambiguously stated.  Finally there is the most exterior and unessential kind of predication, contingent predication.

How can we interpret these notions in modern terms ?

Consider a finitely axiomatisable model $M$ and let $\Gamma(M)$ be the set of all first-order sentences $\phi$ for which $M \Vdash \phi$.  Then $\Gamma(M)$ has a generator $\phi$ such that for any sentence  $\psi \in \Gamma(M)$ we have $\phi \vdash \psi$.  All generators are logically equivalent.  We can divide the sentences in 
$\Gamma(M)$ into generators and non-generators. This might be an modern incarnation of the ancient division between essential and non-essential predication.

We can also consider the formulas which hold of $M$ itself regarded as a set. For instance for $\mathbb{N}$ we have $\forall n. n\in \mathbb{N} \rightarrow n + 1 \in \mathbb{N}$.  Thus $\mathbb{N}$ satisfies the formula $\phi(x)$ given by
$\forall n. n\in x \rightarrow n + 1 \in x$. 
Let the set of these formulas be denoted by $\Delta(M)$.
Suppose that this set is finitely axiomatisable. Then again we have generators and non-generators.  Clearly for a generator formula to be the modern analogue of property or definition or even genus we do not want the constant $M$ to appear in the formula and we want the formula to be structurally irreducible - not equivalent to a subformula of itself.  But how do we decide which formula counts as 'essential' and which does not ? A tentative solution would be that we have \emph{a priori} a certain intuition of the nature of $M$ and that some formulas
have more direct intuitive force than others. Also the definitional nature of a formula must depend on its having the form $ x \in A \&.... $  which reflects the genus + difference template of the ancient theory.

Suppose we had a certain class of mathematical structures $M$ and a finite $G \subset PM$ whose elements ordered by inclusion satisfied the order properties of $\overline{\in}$ of section 2. How do we determine if $G$ should be considered legitimately as coming from the extensions of some system of genera, species and definitions ? This brings us to interesting considerations involving the modern logical notion of (first-order) definability. Note that the set of indefinable real numbers is definable but no individual element of this set can be identified or defined. This is a striking illustration of an infima species and its subordinate individuals.

Up to now we have confronted modern mathematical practice with the ancient theory of definition directly. We now consider this confrontation and comparison as mediated by early modern developments in logic.

With the advent of modernity we acquired a new perspective on definition. Consider the \emph{primitive terms} or \emph{primitive concepts} of Leibniz's \emph{characteristica universalis} project.  These correspond to atomic predicates\footnote{Our model in section 2 also has, as we saw, a notion of 'simple term' defined in terms of indecomposability via $\times$. We also pointed out Aristotle's explicit discussion of complex terms. }. All other concepts are then supposedly defined in terms of these and result from these through logical combination. The idea of 
primitive concepts -  the clear and distinct ideas - was however also anticipated in ancient philosophy\footnote{for instance in the Stoic prolepsis, ennoia and lekta, the \emph{logoi} of Proclus or the 'ineffables' in Socrates' dream in the \emph{Theaetetus }. See \cite{criv} for an interesting discussion on the Stoic theory of definition.}.  The theory of primitive terms departs from the classical theory of definition. Primitive terms are not supreme genera nor infima species. They satisfy \emph{axioms}. Axioms could be contrasted with Aristotelic definitions.   They bridge the gap with modern logic.
We should not speak of a definition in isolation but only of a \emph{system of definitions}. We will return to this in the final section.

We now show how the ancient theory of definition could be reinterpreted in terms of the modern primitive term based approach.

Consider a first-order language with a finite set of contants and $n$-ary predicate symbols.
We then consider a set $\mathcal{C}$ consisting of new constants which are ordered $c_1,...,c_m$ and  a set of new $n$-ary predicate symbols $\mathcal{R}$ which are ordered $R_1,...,R_k$.  Then a \emph{system of definition} consists
of a pair of ordered sets of formulas $(C,R)$ described as follows: $C$ consists
of a sequence $\phi_1,....,\phi_m$ of sentences of the form $c_i = \phi'_i$ where
$\phi'_i$ can only contain constants $c_j$ for $j < i$. $R$
consists in a sequence $\psi_1,...,\psi_k$ of sentences of the form $\forall x_1...x_s. (R_i(x_1,...,x_s) \leftrightarrow \psi'_i)$ where $s$ is the arity of $R_i$ and the only new predicate symbols that may occur in $\psi'_i$ are those in $\psi_j$ for $j < i$. It is natural to impose that all sentences in the definition system be structurally irreducible  as well as satisfying some suitable optimality constraints. For suppose we wanted to define a new predicate symbol $R$ given by a formula $\psi$. Then there must be optimal ways to successively define auxiliary new symbols and build up a definition system ending in the definition of $R$ in the most 'elegant' or 'compact' way possible.
 
Systems of definitions could easily be organised into a tree-like structure of dependency.  We might think of the new symbols as genera or species, immediate species of a genus and the species-genus relation being organized in terms of the \emph{dependency} of the definition system. The \emph{differences} will be the formulas on the right.  And properties will be sentences logically equivalent to differences. Accidents will be formulae which defined constants satisfy but not uniquely, or  formulae $\rho$ such that $\psi'_i \rightarrow \rho$ but not $\rho \rightarrow \psi'_i$ where $\psi'_i$ defines a relation $R$. A striking feature is that the genus-species relation is one of conceptual dependency only, \emph{it has nothing to do with extensionality - a striking departure from the classical theory}\footnote{But consider the case in which we had only monadic predicates $M_i(x)$. Then in prenex normal form a definition will look like $M_i(x) \leftrightarrow Q_1....Q_k x_1...x_k. V_1 \vee ...\vee V_s$ where the $V_j$ are conjunctions of monadic predicates or negations of monadic predicates $M_l$. Then we can obtain an equivalent formula of the form $W_1(x) \vee...\vee W_m(x)$ where
 the $W_i$ are conjunctions of monadic predicates or negations of monadic predicates.}. 
 
 We see thus that the theory of primitive terms might suggest that the ancient  theory of definition while representing a natural and persistent logico-epistemic paradigm has an essentially constructed and relative nature, as an \emph{a posteriori} way of organising (mathematical) knowledge. However consider the following.

 Given a mathematical structure $A$ let $\Gamma(A)$ be the theory of $A$, all valid sentences involving $A$.  Then of course a good candidate for a 'definition' of $A$ is a finite system of axioms $P$ such that all elements in $\Gamma(A)$ can be deduced from $P$. However knowledge of $P$ in itself does not capture the essence of $A$ because there may be important non-obvious truths about $A$, derivable from $P$ which are not immediately evident upon inspecting $P$. They are mediated truths.  Within $\Gamma(A)$ we might identify a finite set of 'remarkable truths'  $E_A$ which can be considered  part of the 'essence' of $A$. And these will likely involve a relation to other structures $B$, $C$, etc. and also belong to $E_B$, $E_C$. Thus the essence of $A$ will be interwoven with the essence of many other mathematical structures. So we have a system of partial essences which can only be (more) fully grasped within the context of a total mathematical theory.
 
 Finally we discuss briefly the practice of definition in type-theoretic foundations for mathematics which are alternatives to the standard set-theoretic foundations and are used in proof assistant software such as Coq and Agda.
 \begin{remark} There is a difference between how definitions are carried out in first-order logic (and we include single and multi-sorted versions) and in systems of higher-order logic. In the first-order case, as we have seen, having standard equality is more convenient for defining constants and functions by means of axioms involving equality. But in the first-order case in general predicate symbols are introduced through axioms involving logical equivalence.  In higher-order systems (in which their is a type for propositions) constants, functions and predicates are treated uniformly and defined all through equality. In (dependent) type theory we have two kinds of equality, an absolute 'external' kind which pertains to term rewriting and an 'internal' kind which is itself a type. Extensionality, the equivalence between internal and external equality has to be postulated. In certain systems this internal equality has a more geometric flavour and resembles the identity between individuals defined in terms of sharing the same infima species. The 'geometry' comes from the fact that we can rely on the intuition that two individuals are 'equal' if they can be continuously deformed into each other. This is the spaces as types paradigm of homotopy type theory.
 	
 \end{remark}
 When formalizing mathematical theories in such software there seem to be two kinds of definition. One kind involves imposing (logical) conditions on a previously defined structure.
 For instance given the concept of category we have the concept of object of a category and thus we can define a terminal object of a category as an object of that category satisfying a certain condition.
  This resembles  classical definitional practice with the genus + difference template. But there is a second kind of definition which is more constructive. We consider (for the sake of simplicity) two basic kinds. One kind involves combining previously defined structures and then imposing a restriction\footnote{we can also think more generally of subtypes and quotient types}.
 For example the concept of natural transformation is certainly not given by a restriction of some previous relevant structure. Rather it is built from a new kind of type\footnote{in Coq it could be, for instance, a 'record'.} involving multiple previously defined structures to which are imposed conditions.  This is related to problems we raised in section 2 about the necessity of considering cartesian products for the category of relations and also the discussion on the problems introduced by the 'comprehension' operation $\times$ on terms. It seems that a missing ingredient in the ancient theory of definition was an explicit theory of the algebraic or categorical structure on genera, they way they can be canonically 'multiplied', combined and merged to form new ones: \emph{necesse est enim genera multiplicanda}. The discussions in\emph{Topics} Book VI, 13-14 (150a-151b) (besides containing interesting insights into Aristotle's mereology) seem to point to a theory of the product of genera as well as quotients. Aristotle distinguishes between simple products such as 'Justice = Temperance $\times$ Courage' and products which take into account how the parts are joined (like for a house). In this case with cannot but think of a quotient type defined on a product of types.
  But let us continue with the second kind of definition which involves recursive definitions and inductive types - both controlled forms of impredicativity. Impredicativity is not alien to ancient philosophy\footnote{For instance there is an impredicative definition of man ('that which is known by everyone') attributed to Democritus discussed by Sextus in the second book of his \emph{Outlines}. See \cite{pro3} for an analysis.}. Inductive types are remarkable in that they are purely intensional rather than extensional definitions of concepts. Intensional in the sense that they tell you how to build elements of a set rather than handing you the entire set all at once.
 Ancient logic was certainly intensional in the sense that it differentiated between concepts having the same extension. For instance, definition and property have the same 'extension' for Aristotle but are entirely different\footnote{It would be interesting to investigate these ancient intensional distinctions in light of the modern problems relating to Frege's puzzle, the fact that in the scope of propositional attitudes we cannot in general switch between a definiendum and its definiens, thus violating Leibniz's law.}.  Frege's puzzle is in fact at the heart of the whole Socratic enterprise. We can 'know' a term and yet not 'know' its definition ( but we do 'know' its property \emph{idion}) . So the scope of the relation of 'knowing' must violate Leibniz's law.

\section{Glimpses of topos theory}

 In this section we end with some suggestions about how the ancient and modern theory might be brought together into a single unifying framework.

The model of section 3 and the discussion at the end of section 4 suggests the following speculations regarding how the ancient theory might be extended and combined with a modern framework. Not only do we need to consider the Cartesian product $X\times X$ of the set $X$ of individual substances (to deal with a kind of relation) but also an arbitrary finite Cartesian product of genera (in programming language terms, in order to define 'records'). This suggests that the following operations are fundamental for a universe of discourse in which definitions can be carried out: \begin{itemize}
	\item Given two sets $X_1, X_2$ to form their product $X_1\times X_2$.
	\item Given a set $Y$ to form the \emph{generalized} power-set $PY$ whose elements are the subsets of $Y$ (not necessarily \emph{all} subsets but including $Y$ itself, the empty set and singletons).
	\item Given a set $Y$ there is a set $\{0,1\}$ such that a certain class of functions $f: Y \rightarrow \{0,1\}$ corresponds bijectively to $PY$.
	\item Given a set $Y$ its elements will also be sets except for the case of a special 'atomic' set $X$ as well as $\{0,1\}$.
\end{itemize}

Given a set $X$ of individuals we could consider the universe $UX$ built up by applying these operations a countable number of times so that $UX$ is the smallest universe closed under the above operations. $UX$ would be a good universe in which clearly nothing resembling Russell's paradox could be derived. In a sense it would be the natural and correct 'theory of sets' which would also be a universe for an extended version of the classical theory of definition. It is interesting to note that if we consider a restricted version of this universe $TX \subset UX$ characterized by containing $X$ and being closed under the operation 

\[ Y_1,...,Y_n \rightarrow P(Y_1\times....\times Y_n)\]

then using the identification of $PY$ with a certain class of functions $f: Y \rightarrow \{0,1\}$ we obtain Russell's simple theory of types based on the type of atoms $X$.

If $TX$ can be defined rigorously, the description of $UX$ is only meant to be a rough heuristic. The great problem was to give a rigorous formal interpretation of these ideas. If the basic notions involved were discovered and explored systematically by Cantor it was - after an unnecessary detour through Russell's paradox and early axiomatic set theory - the profound insight of William Lawvere\footnote{Together with A. Grothendieck in his famous Seminar of Algebraic Geometry (SGA) series.} which succeeded in giving a rigorous category theoretic interpretation of $UX$ and thus of a natural  theory of sets: the concept of \emph{elementary topos}\footnote{For an introduction to topos theory see for instance \cite{Moerdijk}. The present discussion also suggests that it would be interesting to consider the concept of a topos freely generated by a given object $X$.}.  Without going into details we can say very roughly that an elementary topos takes as primitive notions or operations abstract versions of the cartesian product and the powerset operation $P$ and a generalized object $\Omega$ representing $\{0,1\}$ which satisfies the expected correspondence analogous to the third condition above. As suggested by section 3 the repeated applications of $P$ correspond to the progressive ascent to higher genera while the functions $f: P^n X \rightarrow \Omega$ are analogous to differences which carve out the different species. Other constructions in section 3 carry over directly thanks to abstract properties of $P$. On the other hand a topos can also be seen as a model of a theory, as offering a semantics for intuitionistic higher-order logic and thus pertains to contemporary mathematical logic.  We speculate that topos theory might be the place and framework for the ancient and modern approaches to be brought together.

\section{Architecture and genesis of theories}

We have seen how the formalization of the ancient theory of definition presents challenges and problems as well as raising problems with the ancient theory itself and mathematical perspectives regarding this formalization. We then showed how the ancient theory allows us to analyze and ask important questions about the modern practice of definition in mathematics. This approach also carries over to the analysis of the architecture, genesis and epistemology of mathematical theories in general.
We conclude with some considerations about theories which are relevant to the theory of definition. Our goal in this section is to raise questions rather than answer them and to show the connections to Aristotle's \emph{Topics}.
By 'theory' we mean a body of knowledge capable of formal or semi-formal presentation. But we focus on mathematical theories in the framework of axiomatic set theory or type theory. A theory is not merely the set of sentences which can be deduced from a given axiomatic-deductive system. Rather some sentences are marked off as being more meaningful and relevant than others (theorems, lemmas and corollaries and also examples and counterexamples). We call these 'relevant' sentences. Mere logical validities are excluded.  As we saw in section 4, a theory is built up from an dependently ordered (non-circular) hierarchy of definitions. Relevant sentences and the definition system can be organized in a dependency tree which is usually flattened or made linear when the theory is presented\footnote{for simplicity we do not dwell here on the case in which a definition depends on a previous relevant sentence, for instance a uniqueness result}. This intrinsically linear or tree-like structure we call the  \emph{genetic logic} of the theory. How do we justify a genetic logic ? Is it somehow implicit in the axiomatic-deductive system, built into its 'essence' ? Or does the genetic logic\footnote{here genetic logic, unlike for other usages of this term, refers to a purely logical, structural property of the organization of a theory, a particular unfolding and organization of an axiom-deductive system: there is no consideration of historical development.} and the axiomatic-deductive system both derive from some third principle, something that well-known incompleteness results as well as the fact that the 'same' theory can be presented in terms of different axiomatic-deductive systems would naturally suggest ? In what sense are the successive stages of the genetic logic tree not only conceptually dependent on the previous ones but  contained implicitly in them and represent their unfolding ?, or are an unfolding of, the previous elements ? We note sometimes many known definitions turn out to be particular cases of a single more universal definition. For instance the case of Kan extensions in category theory.  

By the \emph{genetic epistemology} of a theory we mean the (optimal) process by which human beings come to learn and understand a theory. This can be either on a individual level or on an historical-cultural level. On an individual level sometimes the genetic logic largely coincides with the genetic epistemology of a theory. Genetic logic can be expressed in the form of research papers, introductions, fundamental treatises and reference works. We wonder to what extent the mirroring of the genetic logic dimension in the epistemological dimension can explain the 'unfolding impulse' mentioned above. How does our knowledge of a concept or theorem (or a particular kind of philosophical reflection on this knowledge)  already in itself lead to the thrust or impulse to find the subsequent concepts or theorems (according to the genetic logic) ? Reflecting on the concept of 'initial object' (or reflecting on our own understanding of such a concept) we could easily be lead to propose the dual concept of 'terminal object' (and vice-versa). Aristotle distinguished between a method starting from things more clear and fundamental 'in themselves' and a method starting from things more clear and fundamental 'to us'. In the \emph{Topics} the condition of being 'better known' is an important ingredient in the theory of definition and property; see Book VI, 4 for a pertinent discussion in which Aristotle shows himself well aware of the relativity of this relation; in 142a he states that what is 'better known' depends on the person and the time. Maybe the concept of Kan extension is more fundamental 'in itself' (relative to category theory) but for cognitive-epistemic reasons 'for us' it is better to first go through a series of concrete cases of this concept. And what is the relationship between historical genetic epistemology and  individual genetic epistemology (both for adults and for the process of child development) ? For instance for centuries children were taught Euclid's \emph{Elements} at school, thus manifesting a kind of cultural law of ontogenetic recapitulation. The genetic logic of the calculus has remained fairly constant since the 17th-century\footnote{Lawvere's synthetic differential geometry wishes to a lead to a even greater conformity to the older way of thinking about infinitesimals.}.
Two additional important questions are: the criteria for the individuality of theories (whether two theories should be considered merely subtheories of a single one) and the theory of analogy between concepts in different theory. The theory of analogy has nothing vague about it. The investigation into analogy is what gave birth to category theory and the notion of functor in the early development of algebraic topology.



By studying a diversity of different theories we can look for conceptual invariants, metatheoretical paradigms which are the basis of analogies between different theories.  In modern mathematics it is category theory that is a shining example of this possibility.

Category theory is itself a theory and its very closely related to logic and type theory. The bare concept of category functions like a supreme genus. As more properties are added these are mirrored in the nature of the internal logic of the category. The way successive relevant properties emerge (and the associated logics become enriched) is certainly not arbitrary but seems to conform some instrinsic meta-theoretic necessity. This is a topic well worth exploring.

\end{document}